\documentclass[12pt,letterpaper,notitlepage]{amsart}
\usepackage{fullpage}

\usepackage{amssymb}
\usepackage{amsmath}
\usepackage{amsthm}
\usepackage{amsfonts}

\usepackage[linktocpage=true]{hyperref} 

\usepackage[numeric,nobysame,lite]{amsrefs}

\usepackage{graphicx}
\usepackage{float}

\usepackage[active]{srcltx} 

\usepackage[applemac]{inputenc}

\newtheorem*{problem}{Problem}
\newtheorem*{solution}{Solution}

\theoremstyle{remark}

\def\Q{\mathbb{Q}}
\def\medio{{\textstyle \frac{1}{2}}\,}
\def\cuarto{{\textstyle \frac{1}{4}}\,}

\begin{document}

\title[Solution to the Gion Shrine Problem]
{A Modern Solution to the Gion Shrine Problem}
\author[Arias de Reyna]{J. Arias de Reyna}
\address{Facultad de Matemáticas \\
Univ.~de Sevilla \\
Apdo.~1160
 \\
41080-Sevilla \\
Spain} 
\thanks{First author supported by  MINECO grant MTM2012-30748.}
\email{arias@us.es}

\author[Clark]{David Clark}
\address{%
Mathematics Department\\
Randolph-Macon College\\
204 Henry Street\\
Ashland, VA 23005\\}
\email{davidclark@rmc.edu}

\date{\today}

\begin{abstract}
We give a new solution to the famous Gion shrine geometry problem from eighteenth-century Japan. Like the classical Japanese solution, ours is given in the form of a degree ten equation. However, our polynomial has the advantage of being much easier to write down. We also provide some additional analysis, including a discussion of existence and uniqueness.
\end{abstract}

\maketitle


\section{Introduction.}

If you had visited Kyoto's Gion shrine\footnote{known today as the Yasaka shrine} around the middle of the eighteenth century, you might have noticed a wooden tablet, inscribed with geometric figures, hanging from one of the eaves. This may not have been a great surprise, since such \emph{sangaku} (literally ``mathematical tablets") were commonplace in temples and shrines throughout Japan at the time. However, this particular tablet happened to hold a problem that would rise to great fame among a generation of Japanese mathematicians.

\begin{problem}
We have a segment of a circle. The line segment $m$ bisects the arc
and chord $AB$.  As shown, we draw a square with side $s$ and an
inscribed circle of  diameter $d$. Let the length $AB = a$. Then,
if
\begin{displaymath}
p = a +
m + s + d\quad\text{ and }\quad q = \frac{m}{a}+ \frac{d}{m}+ \frac{s}{d},
\end{displaymath}
find $a$, $m$, $s$, and $d$ in
terms of $p$ and $q$.
\end{problem}

\begin{figure}[H] \centering
  \includegraphics[width=9cm]{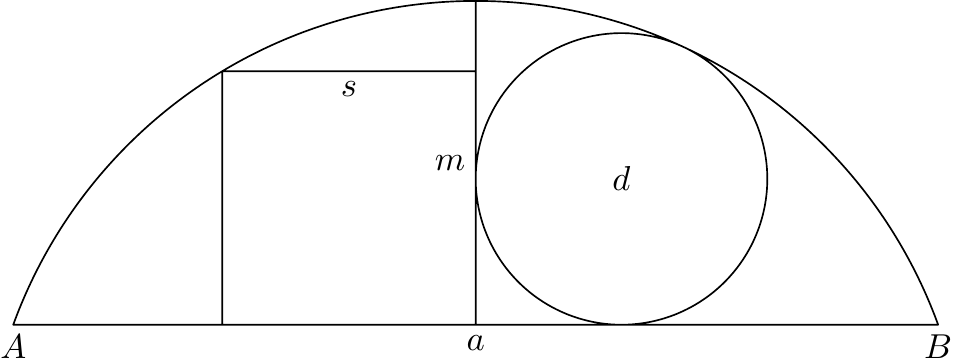}\\
  \caption{Gion shrine problem}\label{figure1}
\end{figure}

\noindent The algebraic particulars of this challenge might strike modern mathematicians as odd, but it has roots in a fascinating niche of mathematical history.

\section{Background}

Eighteenth-century Japan, unified under the Tokugawa shogunate, was a relatively peaceful place where artistic recreation flourished. Those with leisure time indulged in the emerging arts of kabuki and bunraku theater, haiku poetry, \emph{ukiyo-e} woodblock printing, origami---and a homegrown style of mathematics called \emph{wasan}. Practitioners of \emph{wasan} tended to gravitate toward the aesthetics of geometry, and proved wonderful (though esoteric) results about packings of circles, polygons, and ellipses, as well as analogous problems in three dimensions. When a collection of theorems was deemed especially beautiful, it would be inscribed on a \emph{sangaku} and hung in a Buddhist temple or Shinto shrine, both as an offering to the gods and as a challenge to other worshippers.

At the same time, the shogunate's policy of \emph{sakoku} (``closed country") kept Japan intellectually distant from the scientific revolution of the West. The result was an insulated discipline that relied heavily on two sources of established knowledge: the planar geometry results of the Greeks, and the rich body of mathematics imported from China, both of which had long been present in Japanese mathematics. \emph{Sangi} computing rods, a notable Chinese technology, allowed for the numerical computation of roots of polynomials, and were used extensively in Japan. For more about traditional Japanese mathematics, see \cite{MR2423835} and \cite{MR2068248}.

The Gion shrine problem exhibits both geometrical aesthetics and an opportunity to harness the computational power of \emph{sangi}. Tsuda Nobuhisa solved the problem first by deriving a 1024th degree polynomial from whose roots one could derive the result; his solution appeared on a \emph{sangaku} hung from the Gion shrine in 1749. Subsequent progress was made by a mathematician named Nakata, who was able to reduce the necessary polynomial degree to 46. However, a celebrated breakthrough was made by Ajima Naonobu, who in a 1774 handwritten manuscript entitled \emph{Kyoto Gion Gaku Toujyutsu}\footnote{literally ``The Solution to the Kyoto Gion \emph{Sangaku}"} presented a degree ten polynomial solution. Ajima's derivation was first published in 1966 \cite{ajima1}, and has since received a modern analysis that has been translated into English \cite{fukagawa}.

Strikingly, Ajima's approach uses no geometric techniques more sophisticated than the Pythagorean theorem. With a great deal of algebraic persistence, he is able to manipulate a few basic geometric relations into a system of high degree equations in $a$ and $d$. A clever substitution yields four, third degree equations in a single variable $X \neq 0$ whose coefficients are given in terms of $a$, $p$, and $q$. This can be viewed as a homogeneous linear system with nontrivial solution $(X^3, X^2, X, 1)$; any such system must have determinant zero.
Ajima then uses a technique equivalent to Laplace's method of cofactor expansion (c. 1776) to arrive at a degree ten polynomial equation in $a$, which requires nearly a full page to write out completely. It should be noted that, because of \emph{sakoku}, Ajima (1732-1798) may not have even heard of Laplace (1749-1827), and likely was unaware of his results.

The solution given in this paper also has the form of a tenth degree polynomial. In contrast to Ajima's, ours uses the formalism of trigonometry extensively, though it should be noted that these techniques\footnote{and most others present here, excluding the \emph{Mathematica} analysis} would likely have been available to eighteenth century Japanese mathematicians, as well. An upside to this approach is as that it allows for greater geometrical insight---and it results in a polynomial that can comfortably be written in two lines. We also show existence and uniqueness\footnote{Ajima does not seem to have addressed these questions} of solutions, and that, in general, for rational $p$ and $q$, the numbers $a$, $m$, $s$ and $d$ are contained in an extension of $\Q$ of degree 20.

\section{Solution.}

In contrast to Ajima's solution, given in the form of a polynomial in $a$, ours uses a new variable $t$. While $t$ arises somewhat mysteriously from a series of ad-hoc substitutions, ultimately we shall see that in fact $t=d/a$.


\begin{solution}\label{solution}
We start by fixing the constants
\begin{align*}
q_0 & :=-3+\frac{3\sqrt{5}}{2}+\medio\sqrt{\medio(125-41\sqrt{5})}\approx2.3949722\\
\intertext{and}
t_0 & :=\medio\bigl(1-\sqrt{5}+\sqrt{2(5-\sqrt{5})}\bigr) \approx 0.557537.
\end{align*}
Given $p$ and $q$, with $2<q\le q_0$, we first find the unique solution $t\in (0,t_0]$ of the equation
\begin{multline*}\label{equation9}
8t^{10}+(16q-33)t^8+16t^7+(8q^2-49q+56)t^6+(16q-33)t^5-\\-(16q^2-55q+39)t^4
-(16q-22)t^3+(8q^2-23q+18)t^2-t+q-2=0.
\end{multline*}
Then, using $t$, we compute the quantities
\begin{equation*}
\begin{split}
m'&=16t^2\quad d'=16t^2(1-t^2)\quad
a'=16t(1-t^2)\\
s'&=-1+6t^2-t^4+\sqrt{1+20t^2-26t^4+20t^6+t^8}
\end{split}
\end{equation*}
and put $p'=a'+m'+s'+d'$. Finally, the desired quantities will be
\begin{equation*}
a=\frac{p}{p'} a',\quad m=\frac{p}{p'} m',\quad s=\frac{p}{p'}
s',\quad d=\frac{p}{p'} d'.
\end{equation*}
\end{solution}

\begin{figure}[h]
\centering
  \includegraphics[width=8cm]{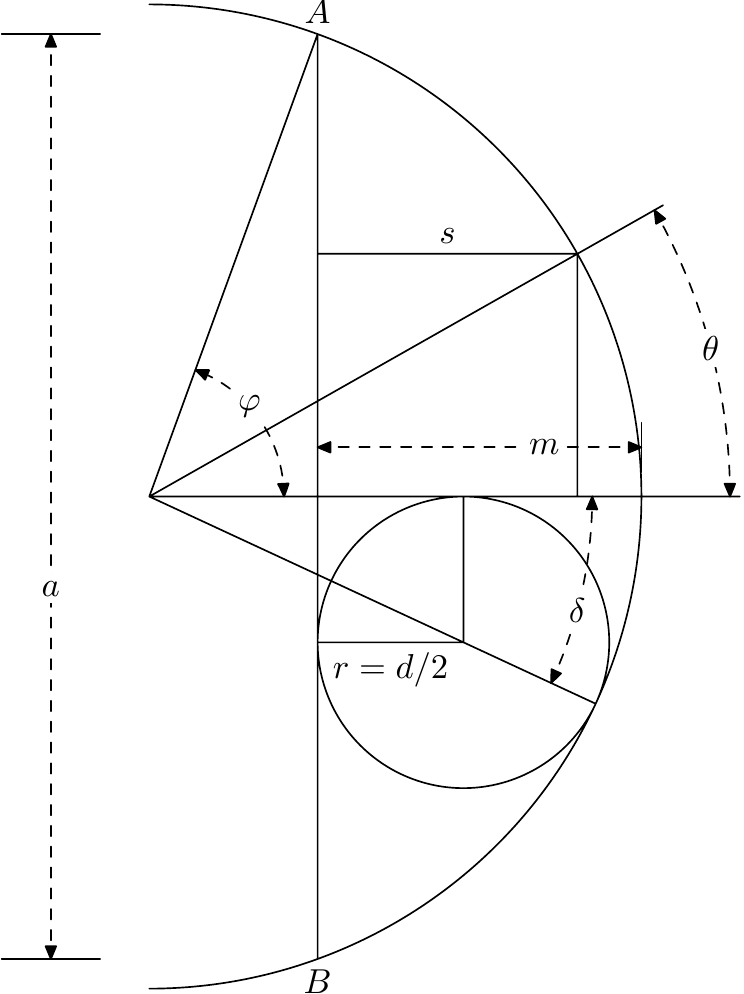}\\
  \caption{Construction for the solution. We first assume that the radius of the circular arc is 1.}
  \label{figure2}
\end{figure}

\begin{proof}

From the definitions of $p$ and $q$, one immediately sees that scaling all lengths by $\lambda$ changes $p$ to $\lambda p$ but leaves $q$ invariant.
As such, the problem is, for all practical purposes, independent of the overall scale. For convenience, let's start by looking for a solution in which the radius of the circular arc is 1.

Observe that the angle $\varphi$ determines the circular segment and, as we shall see, the solution. But for the problem to have a solution, the angle
$\varphi$ is limited  to the interval $0<\varphi\le\varphi_0:=\frac{\pi}{2}+\arctan1/2\approx 117^\circ.$
If $\varphi>\varphi_0$, the square will fail to fit inside the segment. Figure~\ref{figure3} shows the limiting case in which $\varphi=\varphi_0$.
\begin{figure}[h]
\centering
\includegraphics[width=6cm]{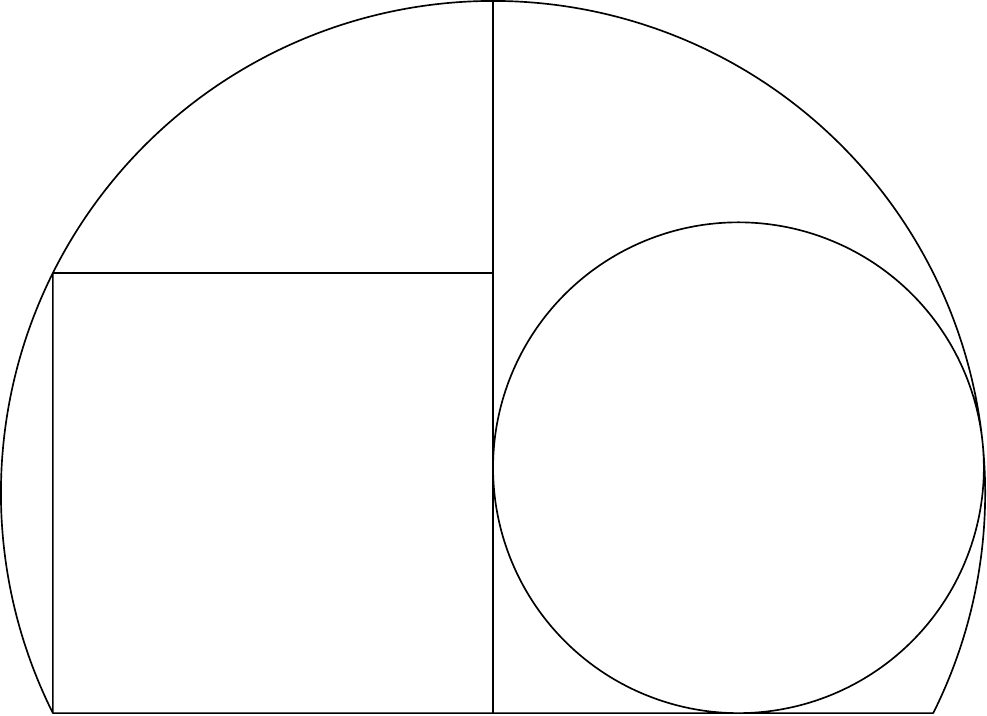}\\
\caption{The extreme case}\label{figure3}
\end{figure}
Not surprisingly, the given value of $q$ completely determines $\varphi$. We shall see this correspondence shortly (Figure~\ref{phi-graph}).

Our first goal will be to express each of the quantities $d$, $m$, $a$, and $s$ in terms of a single variable; the quantity $r=d/2$, the radius of the small circle, happens to be a convenient choice. From Figure~\ref{figure2} we have
\begin{equation}\label{a1}
a = 2\sin\varphi
\end{equation}
\begin{equation}\label{s1}
s= \cos\theta-\cos\varphi=\sin \theta.
\end{equation}

From \eqref{s1}, we find that $\cos\theta=\cos\varphi+\sin\theta$. Thus
\begin{align*}
\cos^2\theta&=\cos^2\varphi+2\cos\varphi\sin\theta+\sin^2\theta,\\
\intertext{and finally}
2\sin^2\theta&+2\cos\varphi\sin\theta-\sin^2\varphi=0.
\end{align*}
This equation, quadratic in $\sin\theta$, has two solutions whose product is
negative. Since in our case $\sin\theta>0$, we must take the greater of the
two solutions,
\begin{displaymath}
\sin\theta=-\medio\cos\varphi+\cuarto\sqrt{4\cos^2\varphi+8\sin^2\varphi},
\end{displaymath}
from which it follows that
\begin{equation}\label{sintheta1}
\sin\theta=\cuarto\sqrt{4+4\sin^2\varphi}-\medio\cos\varphi=
\cuarto\sqrt{8-4\cos^2\varphi}-\medio\cos\varphi.
\end{equation}

The angle $\delta>0$ in Figure~\ref{figure2} can be described by the equations
\begin{equation*}\label{E4}
(1-r)\cos\delta-r=1-m\quad \text{and}\quad \sin\delta=\frac{r}{1-r},
\end{equation*}
where $0<r<1/2$. It follows that
\begin{displaymath}
(1-r)\sqrt{1-\frac{r^2}{(1-r)^2}} - r=1-m,
\end{displaymath}
hence
\begin{align}\label{m1}
m =1+r-\sqrt{1-2r}.
\end{align}
We also have that
\begin{equation*}\label{E6}
m=1-\cos\varphi,
\end{equation*}
from which it follows that
\begin{equation}\label{cosphi}
\cos\varphi=-r+\sqrt{1-2r}.
\end{equation}

Taking a brief digression, we can use \eqref{a1}, \eqref{s1}, \eqref{sintheta1}, and \eqref{cosphi} to write $q$ in terms of $\varphi$. A plot of the implicit function $\varphi(q)$ is given in Figure~\ref{phi-graph}.
\begin{figure}[h]
\centering
\includegraphics[width=10cm]{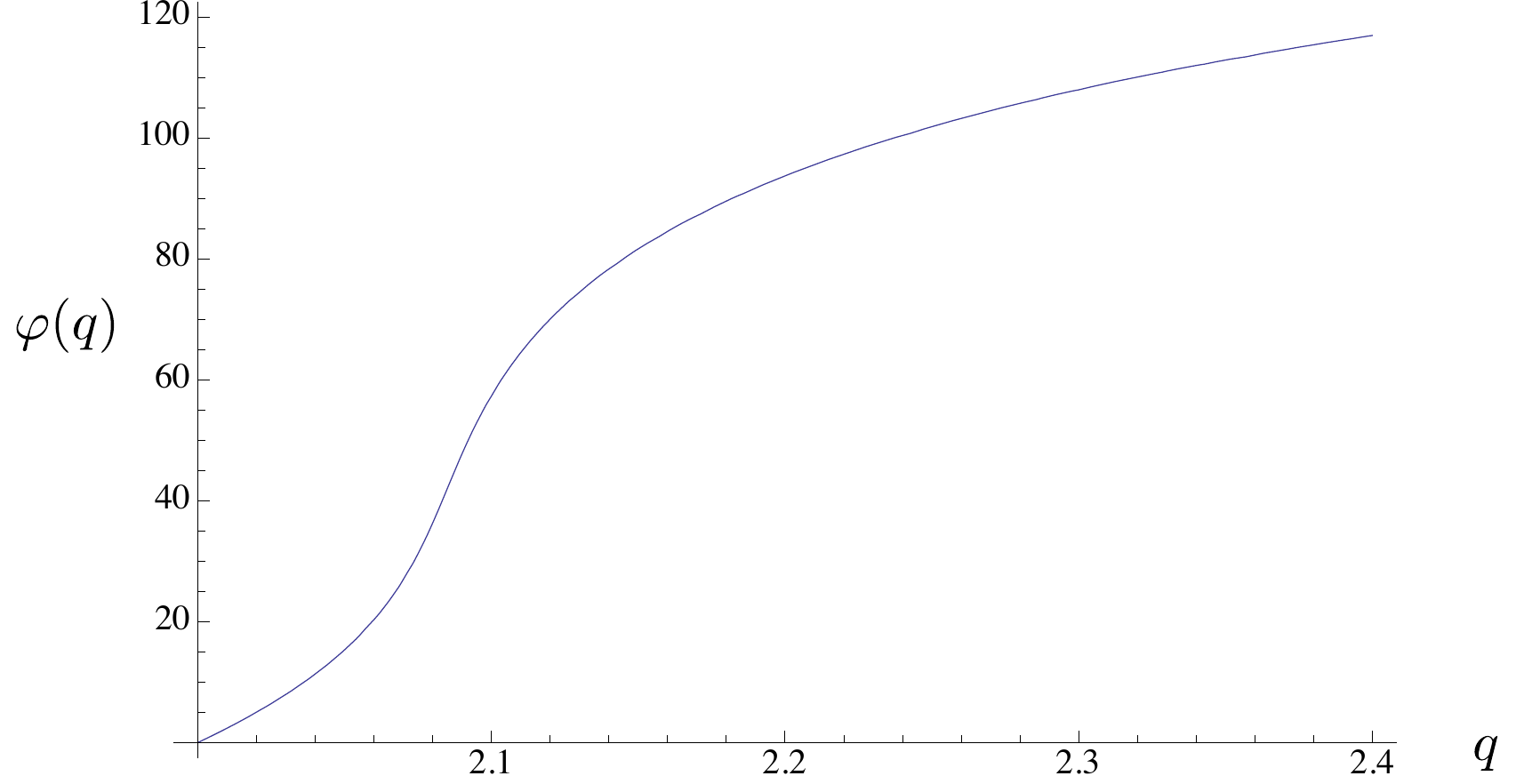}\\
\caption{A given $q$ determines the angle $\varphi$, in degrees.}\label{phi-graph}
\end{figure}

Since $0<\varphi\le\varphi_0$ we have
$1>\cos\varphi\ge-1/\sqrt{5}=\cos\varphi_0$. Then \eqref{cosphi} implies
that $$0<r\le r_0:=-1+\frac{1}{\sqrt{5}}+\sqrt{2-\frac{2}{\sqrt{5}}}.$$
Equation \eqref{cosphi} also tells us that
\begin{equation}
\begin{aligned} \label{sinphi}
\sin^2\varphi&=1-\left((1-2r)+r^2-2r\sqrt{1-2r}\right) \\
&=2r-r^2+2r\sqrt{1-2r},
\end{aligned}
\end{equation}
which, when combined with \eqref{sintheta1} and \eqref{cosphi}, yields
\begin{equation}\label{sintheta2}
\sin\theta=\medio\Bigl(r-\sqrt{1-2r}+\sqrt{1+2r-r^2+2r\sqrt{1-2r}}\Bigr).
\end{equation}

Bringing together our results from \eqref{a1}, \eqref{s1}, \eqref{m1}, \eqref{sinphi}, and \eqref{sintheta2},  we now have the problem's four desired quantities in terms of $r$:
\begin{equation}\label{intermsofr}
\begin{aligned}
d&=2r\\
m&=1+r-\sqrt{1-2r}\\
a&=2\sqrt{2r-r^2+2r\sqrt{1-2r}}\\
s&=\medio\Bigl(r-\sqrt{1-2r}+\sqrt{1+2r-r^2+2r\sqrt{1-2r}}\Bigr)
\end{aligned}
\end{equation}

At this stage, we could write down a polynomial equation relating $q$ and $r$, plug in the given value of $q$ and solve for $r \in (0,r_0]$ (using \emph{sangi}, or \emph{Mathematica}), and use this value to recover the quantities $d$, $m$, $a$, and $s$. However, there are lots of radicals running around in \eqref{intermsofr}, forcing any such polynomial equation to have a large degree. Instead, let's  proceed with a sequence of variable changes that will help to eliminate some of these radicals. First, define $x>0$ by $1-2r=x^2$. The inequality $0 < r \leq r_0 < 1/2$ guarantees that $1>x\ge x_0:=\sqrt{1-2r_0}\approx0.0514622$.  We can write the above
equations \eqref{intermsofr} in terms of $x$:
\begin{equation}\label{intermsofx}
\begin{split}
d&=(1-x^2)\\
m&=\medio(3-2x-x^2)\\
a&=(1+x)\sqrt{3-2x-x^2}\\
s&=\cuarto \left(1-2x-x^2+\sqrt{7+4x-2x^2-4x^3-x^4}\right)
\end{split}
\end{equation}

The situation is better now, but radicals in both $a$ and $s$ would still result in a polynomial degree that's too big. The next natural choice of substitution might be to try a new variable $y^2=3-2x-x^2.$ On its own, such a substitution would create more radicals than it eliminated. However, recognizing that $y^2=3-2x-x^2$ gives the equation of a circle, we can find a rational circle parametrization
in the variable $t$ as below:
\begin{displaymath}
x=\frac{1-3t^2}{1+t^2}, \qquad y=\sqrt{3-2x-x^2}=\frac{4t}{1+t^2}.
\end{displaymath}
This gives us $t^2=(1-x)/(3+x)$, and it is easy to see that when
$x\in[x_0,1)$ we get $t^2\in(0,t_0^{2})$. Thus we'll need to restrict $t$ as follows:
\begin{align}\label{deft0}
0<t\le t_0 :=\medio\bigl(1-\sqrt{5}+\sqrt{2(5-\sqrt{5})}\bigr) \approx 0.557537.
\end{align}

In terms of the new variable $t$ the equations \eqref{intermsofx} become
\begin{equation}\label{E13}
\begin{split}
d&=\frac{8t^2(1-t^2)}{(1+t^2)^2},\quad m=\frac{8t^2}{(1+t^2)^2},\quad a=\frac{8t(1-t^2)}{(1+t^2)^2},\\
s&=\frac{-1+6t^2-t^4+\sqrt{1+20t^2-26t^4+20t^6+t^8}}{2(1+t^2)^2}.
\end{split}
\end{equation}
We can simplify these expressions further by changing the overall scale of the problem: let's choose the radius of the circular arc to be $2(1+t^2)^2$ instead of $1$. In this case, our four quantities can be written
\begin{equation}\label{intermsoft2}
\begin{split}
d&=16t^2(1-t^2),\quad
m=16t^2,\quad
a=16t(1-t^2),\\
s&=-1+6t^2-t^4+\sqrt{1+20t^2-26t^4+20t^6+t^8}.
\end{split}
\end{equation}
Observe that the new variable $t=2r/a$.

Putting everything together,
\begin{align}\label{q(t)}
p&=a+m+s+d \notag \\
&=-1+16t+38t^2-16t^3-17t^4+\sqrt{1+20t^2-26t^4+20t^6+t^8} \notag \\
q&= \frac{m}{a}+ \frac{d}{m}+ \frac{s}{d} \notag \\
&= \frac{-1+22t^2+16t^3-33t^4+16t^6+\sqrt{1+20t^2-26t^4+20t^6+t^8}}{16t^2(1-t^2)}.
\end{align}

Equation \eqref{q(t)}, relating $q$ and $t$, can be rewritten as
\begin{displaymath}
\left(16t^2(-1+t^2)q+(-1+22t^2+16t^3-33t^4+16t^6)\right)^2=(1+20t^2-26t^4+20t^6+t^8),
\end{displaymath}
and further expanded to $32t^2P(t,q)=0$ where
\begin{multline*}
P(t,q)=8t^{10}+(16q-33)t^8+16t^7+(8q^2-49q+56)t^6+(16q-33)t^5-\\-(16q^2-55q+39)t^4
-(16q-22)t^3+(8q^2-23q+18)t^2-t+q-2=0.
\end{multline*}

So, given $p$ and $q$ we first solve for $t$ in $P(t,q)=0$, choosing a root $t\in(0,t_0]$.
Then, equations \eqref{intermsoft2} will recover quantities $a'$, $m'$, $s'$, and $d'$ that correspond to $q$. To get the correct $p$ we need only to rescale the solution by the reciprocal of $p'=a'+m'+s'+d'$. This is precisely the procedure indicated in Solution \ref{solution}.
\end{proof}

\begin{figure}\centering
  \includegraphics[width=8cm]{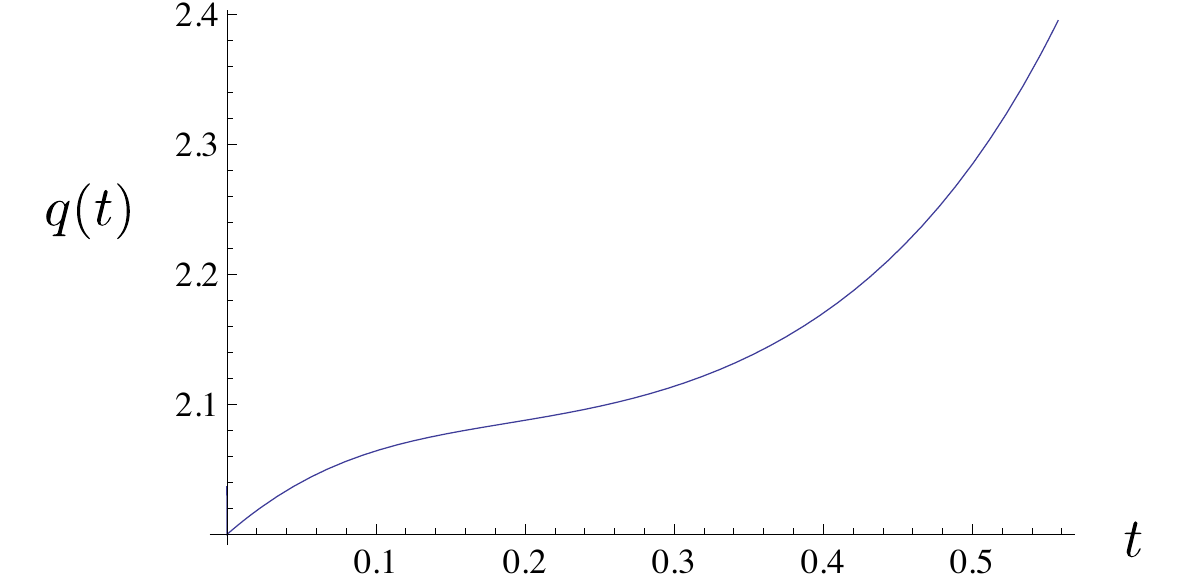}\\
  \caption{Plot of $q(t)$}
  \label{q-graph}
\end{figure}

Not all values of $p$ and $q$ are allowed. We can plot $q$ in terms
of $t$ as given by equation \eqref{q(t)}; see Figure~\ref{q-graph}.
We see that  $q$ varies on the interval $2<q\le q_0$ where
\begin{equation}\label{defq0}
q_0:=q(t_0)=-3+\frac{3\sqrt{5}}{2}+\medio\sqrt{\medio(125-41\sqrt{5})}\approx2.3949722.
\end{equation}
If $q$ is not in the interval $(2, q_0]$, it cannot correspond to an allowable value of $t$; in this case there is no solution to the problem.

On the other hand, $q(t)$ is an increasing function on the interval $[0,t_0]$. To show this,
we compute and rationalize $q'(t)$. The inequality $q'(t)>0$ is equivalent to one of the form $u(t)>0$ for a polynomial $u(t)$; the latter inequality is easy to check. 
Therefore, for each $2<q\le q_0$ there is a unique $0<t \le t_0$ with $P(q,t)=0$.

To summarize, given a pair of positive real numbers $p$ and $q$:
\begin{itemize}
\item if $q\le 2$ or $q>q_0$, there will be no solution to the Gion shrine problem;
\item if $2<q\le q_0$, there will be exactly one solution.
\end{itemize}

For many values of $q\in\Q$, the polynomial $P(q,t)$ is  irreducible. 
For example, this happens  with $q=9/4$.
In such a case, the corresponding $t$ generates a field $\Q(t)$ that is an extension of degree $10$. Here the numbers $a'$, $m'$, $s'$ and $d'$ are contained in $\Q(t)$ or an extension of degree 2 of $\Q(t)$. If we choose $p$ to be rational, the same occurs with the final numbers $a$, $m$, $s$, and $d$.  These numbers could be contained in some subfield. But since $t=2r/a$, this subfield can be only $\Q(t)$ or an extension of degree $2$ of this. Thus, the solutions generate a field extension of degree $10$ (or possibly $20$). 
We conclude with a question: Is there some case in which $a$, $m$, $s$, $d$ are all rational numbers?

\section*{Acknowledgement}
The authors wish to thank Tony Rothman for his interest and helpful comments, and Hidetoshi Fukagawa for sharing a small piece of his life's work on \emph{sangaku} mathematics.


\end{document}